\documentclass[11pt]{article}
\usepackage{amssymb}
\usepackage{amsmath}
\usepackage{theorem}
\usepackage{latexsym}
\usepackage{rotating}
\reversemarginpar
\theoremheaderfont{\bf} 
\theorembodyfont{\sl}
\newtheorem{theo}{Theorem}[section]
{\theorembodyfont{\rm} \newtheorem{defi}[theo]{Definition}}
{\theorembodyfont{\rm} \newtheorem{exa}[theo]{Example}}
{\theorembodyfont{\rm} \newtheorem{rem}[theo]{Remark}}
\newtheorem{prop}[theo]{Proposition}

{\theorembodyfont{\rm}}
{\theorembodyfont{\rm}}
\newenvironment{proof}{{\sc Proof:}}{\mbox{}\hfill$\Box$\par}
\newcommand{\eqnref}[1]{~\mbox{$(${\rm \ref{#1}}$)$}}

\newcommand{\junk}[1]{}
\newcommand{\DS}{\displaystyle}

\newcommand{\N}{{\mathbb N}}
\newcommand{\F}{{\mathbb F}}

\newcommand{\Z}{{\mathbb Z}}
\newcommand{\cC}{{\mathcal C}}

\newcommand{\cS}{{\mathcal S}}

\newcommand{\CC}{convolutional code}
\newcommand{\CCC}{cyclic convolutional code}
\newcommand{\scy}{\mbox{$\sigma$-cyclic}}

\newcommand{\spann}{\mbox{\rm span}\,}
\newcommand{\AutF}{\mbox{${\rm Aut}_{\mathbb F}$}}
\newcommand{\im}{\mbox{\rm im}\,}
\renewcommand{\mod}{\mbox{\rm mod}\,}
\newcommand{\nkdq}{\mbox{$(n,k,\delta)_q$}}
\newcommand{\nkdqm}{\mbox{$(n,k,\delta;m)_q$}}
\newcommand{\nkd}{\mbox{$(n,k,\delta)$}}
\newcommand{\dist}{\mbox{\rm dist}}
\newcommand{\wt}{\mbox{\rm wt}}

\newcommand{\T}{\mbox{$\!^{\sf T}$}}

\newcommand{\ideal}[1]{\mbox{$\langle{#1}\rangle$}}
\newcommand{\Flaurent}{\mbox{$\F(\!(z)\!)$}}
\newcommand{\Azs}{\mbox{$A[z;\sigma]$}}
\newcommand{\Azzs}{\mbox{$A(\!(z;\sigma)\!)$}}
\newcommand{\p}{\mbox{$\mathfrak{p}$}}

\newcommand{\Msigma}{\mbox{${\cal M}^{\sigma}$}}
\newcommand{\plusdot}{\mbox{\small\raisebox{-.8ex}{$\stackrel{+}{\cdot}$}}}

\newcommand{\lideal}[1]{\mbox{$^{^{\bullet\!\!}}\langle{\, #1\, }\rangle$}}

\newcounter{abc}
\newcounter{def}

\newenvironment{alphalist}{\begin{list}{(\alph{abc})\hfill}{\usecounter{abc}
     \topsep-1.4ex \labelwidth.7cm \leftmargin.7cm \labelsep0cm
     \rightmargin0cm \parsep0ex \itemsep.6ex
     \partopsep1.6ex}}{\end{list}}
\newenvironment{arabiclist}{\begin{list}{(\arabic{abc})\hfill}{\usecounter{abc}
     \topsep-1.4ex \labelwidth.7cm \leftmargin.7cm \labelsep0cm
     \rightmargin0cm \parsep0ex \itemsep.6ex
     \partopsep1.6ex}}{\end{list}}

\title{Distance Bounds for Convolutional Codes and Some Optimal Codes}
\date\today
\author{Heide Gluesing-Luerssen$^*$ and Wiland Schmale\footnote{
        Department of Mathematics,
        University of Oldenburg,
        26\,111 Oldenburg,
        Germany, email:
        gluesing@ mathematik.uni-oldenburg.de and
        wiland.schmale@uni-oldenburg.de
        }
        }
\begin{document}
\maketitle

\begin{abstract}
\noindent
After a discussion of the Griesmer and Heller bound for the distance of a
convolutional code we present several codes with various parameters, over various
fields, and meeting the given distance bounds.
Moreover, the Griesmer bound is used for deriving a lower bound for the field size
of an MDS convolutional code and examples are presented showing that, in most
cases, the lower bound is tight.
Most of the examples in this paper are cyclic convolutional codes in a generalized
sense as it has been introduced in the seventies.
A brief introduction to this promising type of cyclicity
is given at the end of the paper in order to make the examples more
transparent.

\end{abstract}


\section{Introduction}
\setcounter{equation}{0}

The fundamental task of coding theory is the construction of good codes, that is,
codes having a large distance and a fast decoding algorithm.
This task applies equally well to block codes and \CC{}s.
Yet, the state of the art is totally different for these two classes of codes.
The mathematical theory of block codes is highly developed and has produced
many sophisticated classes of codes, some of which, like BCH-codes, also come
with an efficient decoding algorithm.
On the other hand, the mathematical theory of \CC{}s is still in the
beginnings.
Engineers make use of these codes since decades, but all \CC{}s used in practice have
been found by systematic computer search and their distances have
been found by computer as well, see for instance~\cite{La73}
and~\cite[Sec.~8]{JoZi99} for codes having the
largest distance among all codes with the same parameters.
Moreover, in all practical situations decoding of \CC{}s is done by
search algorithms, for instance the Viterbi algorithm or one of the sequential decoding
algorithms, e.~g.\ the stack algorithm.
It depends on the algorithm how complex a code may be without exceeding the
range of the decoding algorithms.
However, the important fact about the theory of \CC{}s is that so far no specific
codes are known that allow an {\em algebraic decoding}
(in the present paper a decoding algorithm will be called algebraic if it is capable to
exploit the specific structure of the given code in order to avoid a full search).

Since the seventies quite some effort has been made in order to find algebraic
constructions of \CC{}s that guarantee a large (free)
distance~\cite{Ju73,MCJ73,Ju75,SGR01,GRS03}.
The drawbacks of all these constructions are that, firstly, the field size has to be adapted and
in general becomes quite large and, secondly, so far no algebraic decoding for these codes
is known.
A main feature of most of these constructions is that they make use of
cyclic block codes in order to derive the desired \CC.

Parallel to these considerations there was an independent investigation of
\CC{}s that have a cyclic structure themselves, which also began in the
seventies~\cite{Pi76,Ro79,GSS02,GS02}.
It was the goal of these papers to see whether this additional structure
has, just like for block codes, some benefit for the error-correcting
capability of the code.
The first and very important observation of the seventies was the fact
that a \CC\ which is cyclic in the usual sense is a block
code.
This negative insight has led to a more complex notion of cyclicity for
\CC{}s.
The algebraic analysis of these codes has been completed only
recently in~\cite{GS02} and yields a nice, yet nontrivial, generalization
of the algebraic situation for cyclic block codes.
Furthermore, by now plenty of optimal \CCC{}s have been found in the sense that their
(free) distance reaches the Griesmer bound.
To the best of our knowledge it was, for most cases of the parameters, not known before
whether such optimal codes existed.
Many of these codes are over small fields (like the binary field) and are
therefore well-suited for the existing decoding algorithms.
Along with the algebraic theory of~\cite{GS02} all this indicates that this notion of
cyclicity is not only the appropriate one for \CC{}s but also a very promising one.
Yet, the theory of these codes is still in the beginnings.
So far, no theoretical results concerning the distance of such a code or its
decoding properties are known.
But we are convinced that this class of codes deserves further
investigation and that the theory developed so far will be a good basis
for the next steps.

It is the aim of this paper to present many of these examples in order to
introduce the class of \CCC{}s to the convolutional coding community.
The examples are presented via a generator matrix so that no knowledge about
cyclicity for \CC{}s is required from the reader.
The (free) distances of all these codes have been obtained by a computer program.
A detailed discussion of various distance bounds for convolutional codes over
arbitrary fields shows that all the given codes are optimal
with respect to their distance.
It is beyond the scope of this paper to acquaint the reader with the theory
of \CCC{}s.
However, in Section~\ref{S-CCC} we will give a very brief introduction into
this subject so that the reader may see how the examples have been constructed.
The details of the theory can be found in~\cite{GS02}.

The outline of the paper is as follows.
After reviewing the main notions of convolutional coding theory in the
next section we will discuss in Section~\ref{S-Bounds} various bounds for the free
distance of a \CC, the Griesmer bound, the Heller bound and the
generalized Singleton bound.
The first two bounds are well-known for binary convolutional codes and
can straightforwardly be generalized to codes over arbitrary fields.
It is also shown that for all sets of parameters the Griesmer bound is at
least as good as the Heller bound.
The generalized Singleton bound is an upper bound for the free distance of
a code of given length, dimension, and complexity, but over an arbitrary
field. Just like for block codes a code reaching this bound is called an
MDS code~\cite{RoSm99}.
The Griesmer bound is used for showing how large the field size
has to be in order to allow for an MDS code.
In Section~\ref{S-Exa} many examples of codes are presented reaching the
respective bound.
Most of these examples are \CCC{}s, but we also include some other
codes with the purpose to exhibit certain features of convolutional codes.
For instance, we give examples of MDS codes showing that the lower bounds
for the field size as derived in Section~\ref{S-Bounds} are tight.
Furthermore, an example is given showing that a code reaching the Griesmer
bound may have extreme Forney indices, a phenomenon that does not occur for
MDS codes.
The paper concludes with a brief account of cyclicity for \CC{}s.

\section{Preliminaries}
\setcounter{equation}{0}

We will make use of the following notation.
The symbol~$\F$ stands for any finite field while~$\F_q$ always denotes a field
with~$q$ elements.
The ring of polynomials and the field of formal Laurent series over~$\F$ are
given by
\[
   \F[z]=\Big\{\sum_{j=0}^Nf_jz^j\,\Big|\,N\in\N_0,\,f_j\in\F\Big\}
   \ \text{ and }\
  \Flaurent=\Big\{\sum_{j=l}^{\infty}f_jz^j\,\Big|\,l\in\Z,\,f_j\in\F\Big\}.
\]

The following definition of a convolutional code is standard.

\begin{defi}\label{D-CC}
Let $\F=\F_q$ be a field with~$q$ elements.
An \nkdq-{\em convolutional code\/} is a $k$-dimensional subspace~$\cC$ of the
vector space $\Flaurent^n$ of the form
\[
  \cC=\im G:=\big\{uG\,\big|\, u\in\Flaurent^k\big\}
\]
where $G\in\F[z]^{k\times n}$ satisfies
\begin{alphalist}
\item $G$ is {\sl right invertible}, i.~e.
      there exists some matrix $\tilde{G}\in\F[z]^{n\times k}$ such that
      $G\tilde{G}=I_k$.
\item $\delta=\max\{\deg\gamma\mid \gamma$ is a $k$-minor of $G\}$.
\end{alphalist}
We call $G$ a {\sl generator matrix\/} and~$\delta$ the {\sl complexity\/}
of the code~$\cC$.
\end{defi}

The complexity is also known as the {\em overall constraint length\/}~\cite[p.~55]{JoZi99}
or the {\em degree\/}~\cite[Def.~3.5]{McE98} of the code.
Notice that a generator matrix is always polynomial and has a
polynomial right inverse.
This implies that in the situation of Definition~\ref{D-CC} the polynomial
codewords belong to polynomial messages, i.~e.
\begin{equation}\label{e-cpoly}
  \cC\cap\F[z]^n=\big\{uG\,\big|\,u\in\F[z]^k\big\}.
\end{equation}
In other words, the generator matrix is delay-free and non-catastrophic.
As a consequence, a convolutional code is always uniquely determined by its polynomial
part.
Precisely, if $\cC=\im G$ and $\cC'=\im G'$ where
$G,\,G'\in\F[z]^{k\times n}$ are right invertible, then
\begin{equation}\label{e-cpolyunique}
  \cC=\cC'\Longleftrightarrow \cC\cap\F[z]^n=\cC'\cap\F[z]^n.
\end{equation}
This follows from\eqnref{e-cpoly} and the fact that
$\{uG\mid u\in\F[z]^k\}=\{uG'\mid u\in\F[z]^k\}$ is equivalent to
$G'=VG$ for some matrix $V\in\F[z]^{k\times k}$ that is invertible over
$\F[z]$.
This also shows that the complexity of a code
does not depend on the choice of the generator matrix.
From all this it should have become clear that with respect to code
construction there is no difference whether one works in the context of
infinite message and codeword sequences (Laurent series) or finite ones
(polynomials) as long as one considers right invertible generator
matrices.
Only for decoding it becomes important whether or not one may assume the
sent codeword to be finite.
The issue whether convolutional coding theory should be based
on finite or infinite message sequences, has first been raised and discussed in
detail in~\cite{RSY96,Ro01}.

It is well-known~\cite[Thm.~5]{Fo70} or~\cite[p.~495]{Fo75} that each convolutional code
has a minimal generator matrix in the sense of the next definition.
In the same paper~\cite[Sec.~4]{Fo75} it has been shown how to derive such a matrix
from a given generator matrix in a constructive way.

\begin{defi}
\begin{arabiclist}
\item For $v=\sum_{j=0}^N v_jz^j\in\F[z]^n$ where $v_j\in\F^n$ and $v_N\not=0$
      let $\deg v:=N$ be the {\em degree\/} of~$v$.
      Moreover, put $\deg0=-\infty$.
\item Let $G\in\F[z]^{k\times n}$ be a right invertible matrix with complexity
      $\delta=\max\{\deg\gamma\mid \gamma$ is a $k$-mi\-nor of $G\}$
      and let $\nu_1,\ldots,\nu_k$ be the degrees of the rows of~$G$ in the sense of~(1).
      We say that $G$ is {\sl minimal\/} if $\delta=\sum_{i=1}^k\nu_i$.
      In this case, the row degrees of~$G$ are uniquely determined by the code
      $\cC:=\im G\subseteq\Flaurent^n$. They are called the {\sl Forney indices\/}
      of~$\cC$ and the number $\max\{\nu_1,\ldots,\nu_k\}$ is said to be the {\sl memory\/}
      of the code. An \nkdq-code with memory~$m$ is also called an
      \nkdqm-code.
\end{arabiclist}
\end{defi}

From the above it follows that an \nkdq-convolutional code has a constant generator matrix if
and only if $\delta=0$. In that case the code can be regarded as an $(n,k)_q$-block
code.

The definition of the distance of a convolutional code is straightforward.
For a constant vector $w=(w_1,\ldots,w_n)\in\F^n$ we define its {\em (Hamming) weight\/} as
$\wt(w)=\#\{i\mid w_i\not=0\}$.
For a polynomial vector $v=\sum_{j=0}^N v_j z^j\in\F[z]^n$, where $v_j\in\F^n$,
the {\em weight\/} is defined as $\wt(v)=\sum_{j=0}^N\wt(v_j)$.
Then the {\em (free) distance\/} of a code $\cC\subseteq\Flaurent^n$ with generator matrix
$G\in\F[z]^{k\times n}$ is given as
\[
   \dist(\cC):=\min\big\{\wt(v)\,\big|\, v\in\cC\cap\F[z]^n,\;v\not=0\big\}.
\]
By virtue of\eqnref{e-cpoly} this can be rephrased as
$\dist(\cC)=\min\{\wt(uG)\mid u\in\F[z]^k,\,u\not=0\}$.

When presenting some optimal codes in Section~\ref{S-Exa} we will also investigate the column
distances of the codes.
For each $l\in\N_0$ the $l$th {\em column distance\/} of~$\cC$ is defined as
\begin{equation}\label{e-coldist}
  d^c_l=\min\Big\{\wt\big((uG)_{[0,l]}\big)\,\Big|\,
                      u\in\F[z]^k,u_0\not=0\Big\}
\end{equation}
where for a polynomial vector $v=\sum_{j=0}^N v_jz^j$ we define
$v_{[0,l]}=\sum_{j=0}^{\min\{N,l\}} v_jz^j$.
It can easily be shown~\cite[Thm.~3.4]{JoZi99} that for each code~$\cC$ there exists some $M\in\N_0$
such that
\begin{equation}\label{e-coldistfree}
    d^c_0\leq d^c_1\leq d^c_2\ldots\leq
    d^c_M=d^c_{M+1}=\ldots=\dist(\cC).
\end{equation}

\section{Distance Bounds}\label{S-Bounds}
\setcounter{equation}{0}

In this section we want to present some upper bounds for the distance of a convolutional
code.
These bounds are quite standard for binary convolutional codes and can be found in
Chapter~3.5 of the book~\cite{JoZi99}.
The proof for arbitrary fields goes along the same lines of arguments, but for sake
of completeness we wish to repeat the arguments in this paper.
We will also compare the numerical values of the bounds with each other.

Let us begin with recalling various distance bounds for block codes.
The Plotkin bound as given below can be found in~\cite[1.4.3]{Be98}, but
can also easily be derived from the more familiar formula
\begin{equation}\label{e-Plot}
   \text{ if }d>\theta n\text{ where }\theta=\frac{q-1}{q},
   \text{ then }
   q^k\leq\frac{d}{d-\theta n},
\end{equation}
see for instance~\cite[(5.2.4)]{vLi99}.
As for the Singleton and the Griesmer bound we also refer
to~\cite[Ch.~5.2]{vLi99}.

\begin{theo}\label{T-BCB}
Let $\cC\subseteq\F^n$ be an $(n,k)_q$-block code and let $d=\dist(\cC)$.
Then
\[
\begin{array}{ll}
    {\DS d\leq n-k+1} &\quad \text{(Singleton bound),}\\[1.7ex]
    {\DS d\leq\Big\lfloor\frac{nq^{k-1}(q-1)}{q^k-1}\Big\rfloor}&\quad
       \text{(Plotkin bound),}\\[1.7ex]
    {\DS\sum_{l=0}^{k-1}\Big\lceil\frac{d}{q^l}\Big\rceil\leq n}&\quad
        \text{(Griesmer bound).}
\end{array}
\]
An $(n,k)_q$-code~$\cC$ with $\dist(\cC)=n-k+1$ is called an MDS code.
\end{theo}

Notice that the Singleton bound does not take the field size into account.
As a consequence the question arises as to how large the field size~$q$ has to be in
order to allow the existence of MDS codes and how to construct such codes.
Answers in this direction can be found in~\cite[Ch.~11]{MS77}.

It is certainly well-known that the Griesmer bound is at least as good as the
Plotkin bound.
The importance of the Plotkin bound, however, is that it also applies to nonlinear
block codes, in which case it is usually given as in\eqnref{e-Plot} and
with $M:=|\cC|$ instead of $q^k$.
Since we did not find a comparison of the two bounds for linear block codes
in the literature we wish to present a short proof of this statement.
We also include the relation between the Griesmer and the Singleton bound.

\begin{prop}\label{P-GP}
Given the parameters $n,\,k,\,d$, and $q\in\N$ where $k<n$ and $q$ is a prime
power. Assume $\sum_{l=0}^{k-1}\Big\lceil\frac{d}{q^l}\Big\rceil\leq n$.
Then
\begin{alphalist}
\item ${\DS  d\leq\Big\lfloor\frac{nq^{k-1}(q-1)}{q^k-1}\Big\rfloor}$,
\item $d\leq n-k+1$.
\end{alphalist}
\end{prop}
There is no relation between the Plotkin and the Singleton bound in this
generality.
Roughly speaking, for relatively large values of~$q$ the Singleton bound is
better than the Plotkin bound while for small values the Plotkin bound is
better.

\begin{proof}
(a) Assume to the contrary that
$d>\big\lfloor\frac{nq^{k-1}(q-1)}{q^k-1}\big\rfloor$. Since~$d$ is an
integer this implies that $d>\frac{nq^{k-1}(q-1)}{q^k-1}$.
Thus
\[
  \sum_{l=0}^{k-1}\Big\lceil\frac{d}{q^l}\Big\rceil
  \geq\sum_{l=0}^{k-1}\frac{d}{q^l}
  >\sum_{l=0}^{k-1}\frac{n(q-1)}{q^k-1}q^{k-1-l}
  =\frac{n(q-1)}{q^k-1}\sum_{l=0}^{k-1}q^l=n.
\]
(b) follows from
$\sum_{l=0}^{k-1}\big\lceil\frac{d}{q^l}\big\rceil\geq d+k-1$.
\end{proof}

One should also recall that the Griesmer bound is not tight.
An example is given by the parameters $n=13,\,k=6,\,q=2$ in which case the
Griesmer bound shows that the distance is upper bounded by~$5$.
But it is known that no $(13,6)_2$-code with distance~$5$ exists,
see~\cite[p.~69]{vLi99}.

We will now present the generalization of these bounds to convolutional codes.
Let us begin with the Singleton bound.
The following result has been proven in~\cite[Thm.~2.2]{RoSm99}.

\begin{theo}\label{T-MDSC}
Let $\cC\subseteq\Flaurent^n$ be an \nkd-code. Then
\begin{alphalist}
\item The distance of~$\cC$ satisfies
      \[
         \dist(\cC)\leq(n-k)\Big(\Big\lfloor\frac{\delta}{k}\Big\rfloor+1\Big)
               +\delta+1=:S\nkd.
      \]
      The number $S\nkd$ is called the generalized Singleton bound for the
      parameters \nkd\ and we call the code~$\cC$ an MDS code if
      $\dist(\cC)=S\nkd$.
\item If~$\cC$ is an MDS code and $\delta=ak+r$ where $a\in\N_0$ and $0\leq r\leq k-1$,
      then the Forney indices of~$\cC$ are given by
      \[
          \underbrace{a,\ldots,a}_{k-r\text{ times}},
          \underbrace{a+1,\ldots,a+1}_{r\text{ times}}.
      \]
      Hence the code is compact in the sense of~\cite[Cor.~4.3]{McE98}.
\end{alphalist}
\end{theo}

Just like for block codes the acronym MDS stands for maximum distance
separable.
In~\cite[Thm.~2.10]{RoSm99} it has been shown that for all given parameters
$n,\,k,\,\delta$ and all primes~$p$ there exists an MDS code over a suitably large field
of characteristic~$p$.
The proof is non-constructive and, as a consequence, does not give a hint about the field
size required.
In~\cite[Thm.~3.3]{SGR01} a construction of \nkd-MDS codes over fields~$\F_{p^r}$ is
given under the condition that
$n|(p^r-1)$ and $p^r\geq\frac{n\delta^2}{k(n-k)}$.
Notice that this requires~$n$ and the characteristic~$p$ being coprime.
This result gives first information about the field size required in order to guarantee
the existence of an MDS code.
However, many examples of MDS codes over smaller fields are known.
We will present some of them in the next section.
Although they all have a certain structure in common (they are cyclic in the sense of
Section~\ref{S-CCC}) we do not know any general construction for cyclic MDS codes yet.

Now we proceed with a generalization of the Plotkin and Griesmer bound to
convolutional codes.

\begin{theo}\label{T-CCB}
Let $\cC$ be an \nkdqm-convolutional code having distance $\dist(\cC)=d$.
Moreover, let
\[
   \hat{\N}=\left\{\begin{array}{ll}\N:=\{1,2,\ldots\},&\text{if }km=\delta\\[.6ex]
                          \N_0:=\{0,1,2,\ldots\},&\text{if }km>\delta\end{array}\right.
\]
Then
$$
\begin{array}{ll}
   d&\leq{\DS \min_{i\in\hat{\N}}
         \Big\lfloor\frac{n(m+i)q^{k(m+i)-\delta-1}(q-1)}{q^{k(m+i)-\delta}-1}\Big\rfloor}=:H_q(n,k,\delta;m)
     \hfill \text{(Heller bound)}\\[3ex]
   d&\leq{\DS\max\Big\{d'\in\{1,\ldots,S\nkd\}\,\Big|\,\sum_{l=0}^{k(m+i)-\delta-1}
                    \Big\lceil\frac{d'}{q^l}\Big\rceil\leq n(m+i)
         \text{ for all }i\in\hat{\N}\Big\}}\\
    &=:G_q(n,k,\delta;m)\hfill \text{(Griesmer bound)}
\end{array}
$$
Moreover, $G_q(n,k,\delta;m)\leq H_q(n,k,\delta;m)$.
\end{theo}

In the binary case ($q=2$) both bounds can be found
in~\cite[3.17~and~3.22]{JoZi99}.
In that version the first bound has been proven first by Heller
in~\cite{He68}.
The Griesmer bound as given above differs slightly from the one given
at~\cite[3.22]{JoZi99}.
We have upper bounded the possible values for~$d'$ by the generalized Singleton
bound, which is certainly reasonable to do.
As a consequence, the Griesmer bound is always less than or equal to the generalized
Singleton bound.
This would not have been the case had we taken the maximum over all $d'\in\N$.
This can be seen by taking the parameters $\nkdqm=(5,2,3;3)_8$. In this case
the generalized Singleton bound is $S\nkd=10$ but the inequalities of the Griesmer
bound are all satisfied for the value $d'=12$.

The proof of the inequalities above is based on the same idea as in the
binary case as we will show now.

\begin{proof}
The last statement follows from Proposition~\ref{P-GP}(a).
As for the bounds themselves we will see that they are based on certain block codes
which appear as subsets of the given
convolutional code~$\cC$. This will make it possible to apply the block code bounds
of Theorem~\ref{T-BCB}.
The subcodes to be considered are simply the subsets of all codewords
corresponding to polynomial messages with an upper bounded degree.
\\
Let $\cC=\im G$, where $G\in\F[z]^{k\times n}$ is right-invertible and minimal with
Forney indices $\nu_1,\ldots,\nu_k$. Hence $\delta=\sum_{i=1}^k \nu_i$ and
$m=\max\{\nu_1,\ldots,\nu_k\}$.
Notice that $km\geq\delta$ and $km=\delta\Longleftrightarrow \nu_1=\ldots=\nu_k=m$.
For each $i\in\N_0$ define
\[
  U_i=\{(u_1,\ldots,u_k)\in\F[z]^k\mid \deg u_l\leq m+i-1-\nu_l\text{ for }
  l=1,\ldots,k\}.
\]
This implies $u_l=0$ if $\nu_l=m$ and $i=0$.
In particular, $U_i=\{0\}\Longleftrightarrow km=\delta\text{ and }i=0$ and
this shows that $i=0$ has to be excluded if $km=\delta$.
Obviously, the set $U_i$ is an $\F$-vector space and
$\dim_{{\mathbb F}} U_i=\sum_{l=1}^k(m+i-\nu_l)=k(m+i)-\delta$.
Consider now $\cC_i:=\{uG\mid u\in U_i\}$ for $i\in\N_0$.
Then $\cC_i\subseteq\cC$ and~$\cC_i$ is an $\F$-vector space and, by injectivity
of~$G$,
\[
   \dim_{\mathbb F}\cC_i=\dim_{{\mathbb F}} U_i=k(m+i)-\delta.
\]
Furthermore, minimality of the generator matrix~$G$ tells us that
\[
  \deg(uG)=\max_{l=1,\ldots,k}(\deg u_l+\nu_l)\leq m+i-1
  \text{ for all }u\in U_i,
\]
see~\cite[p.~495]{Fo75}.
Hence $\cC_i$ can be regarded as a block code of length~$n(m+i)$ and dimension
$k(m+i)-\delta$ for all $i\in\hat{\N}$.
Since $\dist(\cC)\leq\dist(\cC_i)$ for all $i\in\hat{\N}$ we obtain the desired
results by applying the Plotkin and Griesmer bounds of Theorem~\ref{T-BCB} to the
codes~$\cC_i$.
\end{proof}

The proof shows that the existence of an \nkdqm-code meeting the Griesmer bound
implies the existence of $(n(m+i),k(m+i)-\delta)_q$-block codes having at least the same distance
for all $i\in\hat{\N}$.
The converse, however, is not true, since the block codes have to have some
additional structure.
We will come back to this at the end of this section.

One should note that these bounds do only take the largest Forney index, the memory,
into account. More precisely, the proof shows that codewords having degree smaller
than~$m-1$ are never taken into consideration.
As a consequence, codes with a rather bad distribution of the Forney indices will
never attain the bound. For instance, for a code with parameters
$(n,k,\delta;m)_q=(5,3,4;2)_2$ the Griesmer bound shows that the
distance is upper bounded by~$6$.
This can certainly never be attained if the Forney indices of that code
are given by $0,2,2$ since in that case a constant codeword exists. Hence
the Forney indices have to be $1,1,2$. In this case a code with
distance~$6$ does indeed exist, see the first code given in Table~I of
Section~\ref{S-Exa}.
But also note that, on the other hand, a code reaching the Griesmer bound need not
be compact (see Theorem~\ref{T-MDSC}(b)); an example is given by the
$(5,2,6;4)_2$-code given in Table~I of the next section.

The Griesmer bound as given above has the disadvantage that infinitely many
inequalities have to be considered.
A simple way to reduce this to finitely many inequalities is obtained by making
use of the generalized Singleton bound $S\nkd$.
Instead of this bound one could equally well use any of the numbers occurring
on the right hand side of the Heller bound.

\begin{prop}\label{P-Gfinite}
Given the parameters $n,\,k,\,m,\,\delta$ such that $k<n$ and
$km\geq\delta$ and let~$q$ be any prime power.
Define the set~$\hat{\N}$ as in Theorem~\ref{T-CCB}.
Furthermore, let $i_0\in\N$ be such that $q^{k(m+i_0)-\delta}\geq S\nkd$ and
put $\hat{\N}_{\leq i_0}:=\hat{\N}\cap\{0,1,\ldots,i_0\}$.
Then
\begin{equation} \label{e-Gfinite}
  G_q(n,k,\delta;m)\!=\!\max\Big\{d'\in\{1,\ldots, S\nkd\}\,\Big|\!
     \sum_{l=0}^{k(m+i)-\delta-1}\!\Big\lceil\frac{d'}{q^l}\Big\rceil
     \leq n(m+i)\text{ for all }i\in\hat{\N}_{\leq i_0}\Big\}.
\end{equation}
Hence the distance of an $\nkdqm$-code is upper bounded by the number given
in\eqnref{e-Gfinite}.
\end{prop}

We will see in the next section that the Griesmer bound is tight for many sets of
parameters.

{\sc Proof:}
Notice that for $a\geq S\nkd$ we have
$\big\lceil\frac{d'}{a}\big\rceil=1$ since $d'\leq S\nkd$.
As for\eqnref{e-Gfinite} it suffices to show that whenever~$d'$ satisfies the
inequality
$\sum_{l=0}^{k(m+i)-\delta-1}\big\lceil\frac{d'}{q^l}\big\rceil\leq n(m+i)$
for some $i\geq i_0$, then it also satisfies the inequality for $i+1$.
But this follows easily from
$$
 \sum_{l=0}^{k(m+i+1)-\delta-1}\Big\lceil\frac{d'}{q^l}\Big\rceil
 =\sum_{l=0}^{k(m+i)-\delta-1}\Big\lceil\frac{d'}{q^l}\Big\rceil
 +\sum_{l=k(m+i)-\delta}^{k(m+i+1)\delta-1}\Big\lceil\frac{d'}{q^l}\Big\rceil
 \leq n(m+i)+k\leq n(m+i+1).
 \eqno\Box
$$

The finite sets for~$d'$ and~$i$ in\eqnref{e-Gfinite} are not optimized, but
they are good enough for our purposes since they allow for a computation of the
Griesmer bound in finitely many steps.
Unfortunately,\eqnref{e-Gfinite} does not reveal the block code case where
only the index $i=1$ has to be considered according to Theorem~\ref{T-BCB}.
The consistency of the Griesmer bound for $m=\delta=0$ with that
case is guaranteed by the following result.

\begin{prop}\label{P-Gblock}
Given the parameters $n,\,k$, and~$q$.
Then
\[
  \max\Big\{d'\in\N\,
    \Big|\,\sum_{l=0}^{ki-1}\Big\lceil\frac{d'}{q^l}\Big\rceil
    \leq ni\text{ for all }i\in\N\Big\}
  =\max\Big\{d'\in\N\,\Big|\,
  \sum_{l=0}^{k-1}\Big\lceil\frac{d'}{q^l}\Big\rceil\leq n\Big\}.
\]
\end{prop}

\begin{proof}
Let~$d'$ be any number satisfying
$\sum_{l=0}^{k-1}\big\lceil\frac{d'}{q^l}\big\rceil\leq n$.
We have to show that~$d'$ satisfies the inequalities given on the
left hand side for all $i\in\N$.
In order to do so, notice that according to Proposition~\ref{P-GP}(a)
\[
  d'\leq\Big\lfloor\frac{nq^{k-1}(q-1)}{q^k-1}\Big\rfloor
  \leq\frac{nq^{k-1}}{1+q+\ldots+q^{k-1}}\leq\frac{n}{k}q^{k-1}.
\]
But this implies $\frac{d'}{q^l}<\frac{n}{k}$ for all
$l\geq k$, thus
$\big\lceil\frac{d'}{q^l}\big\rceil\leq\frac{n}{k}$ and
\[
  \sum_{l=0}^{ki-1}\Big\lceil\frac{d'}{q^l}\Big\rceil
  =\sum_{l=0}^{k-1}\Big\lceil\frac{d'}{q^l}\Big\rceil
   +\sum_{l=k}^{ki-1}\Big\lceil\frac{d'}{q^l}\Big\rceil
  \leq n +k(i-1)\frac{n}{k}
  =ni.
\]
This proves the assertion.
\end{proof}


Finally we want to investigate as to how big the field size~$q$ has to be
in order to allow for an MDS code with parameters $\nkdq$.
A first estimate can be achieved by using the Griesmer bound in
combination with the generalized Singleton bound.

\begin{theo}\label{T-MDSfieldsize}
Let $\cC\subseteq\Flaurent^n$ be an $\nkdqm$-MDS code, thus
$d:=\dist(\cC)=S\nkd=(n-k)\big(\big\lfloor\frac{\delta}{k}\big\rfloor+1\big)+\delta+1$.
Then the field size~$q$ satisfies
\[
    q\geq\left\{\begin{array}{cl}
            \frac{d}{n-k+1},&\text{if }[k=1] \text{ or }[k>1\text{ and }km=\delta+1]
            \\[1ex]
            d,&\text{if }[k>1 \text{ and }km\not=\delta+1].
           \end{array}\right.
\]
\end{theo}
The estimate above also covers the block code case as given
in~\cite[p.~321]{MS77}.

\begin{proof}
We will consider the various cases separately.
In each case we will apply the inequality
\begin{equation}\label{e-GriesMDS}
  \frac{d}{q}\leq
  n(m+i)-d-\sum_{l=2}^{k(m+i)-\delta-1}\Big\lceil\frac{d}{q^l}\Big\rceil,
\end{equation}
which is a simple consequence of the Griesmer bound, to the case $d=S\nkd$.
Moreover we will make use of the fact that $\big\lceil\frac{d}{q^l}\big\rceil\geq1$
for all $l\in\N$.

\underline{$k=1$:}\quad
In this case $m=\delta$ and $d=n(m+1)$.
Since $k(m+i)-\delta-1=i-1$ Inequality\eqnref{e-GriesMDS} gives us
\[
  \frac{d}{q}
  \leq n(m+i)-n(m+1)-(i-2)
  =n(i-1)-i+2
\]
for all $i\geq2$. This shows $q\geq\frac{d}{n}$ as desired.
Using $i=1$ in the Griesmer bound simply leads to $d\leq n(m+1)$.
This is true by assumption and gives no further condition on~$q$.
\\
\underline{$k>1$ and $km=\delta$:}\quad
Now $m=\frac{\delta}{k}$ and thus $d=(n-k)(m+1)+mk+1$.
Using $k(m+i)-\delta-1=ki-1$ we obtain from Inequality\eqnref{e-GriesMDS}
\[
  \frac{d}{q}
  \leq n(m+i)-(n-k)(m+1)-mk-1-(ki-2)
  =(n-k)(i-1)+1
\]
for all $i\geq1$.
Using $i=1$ leads to $q\geq d$.
\\
\underline{$k>1$ and $km>\delta$:}\quad
In this case $m=\big\lfloor\frac{\delta}{k}\big\rfloor+1$, see Theorem~\ref{T-MDSC}(b),
and $d=(n-k)m+\delta+1$.
Therefore Inequality\eqnref{e-GriesMDS} leads to
\[
  \frac{d}{q}\leq
  n(m+i)-(n-k)m-\delta-1-\big(k(m+i)-\delta-2\big)
  =(n-k)i+1
\]
for all $i\geq1$.
This shows $q\geq\frac{d}{n-k+1}$.
In order to finish the proof we have to consider also $i=0$.
In the case $km=\delta+1$ the Griesmer bound applied to $i=0$
simply leads to $d\leq nm$, which is true anyway, and no additional
condition on~$q$ arises.
If $km-\delta>1$ a better bound can be achieved.
Since $\big\lfloor\frac{\delta}{k}\big\rfloor=m-1$, we obtain after division
with remainder of $\delta$ by~$k$ an identity of the form $\delta=(m-1)k+r$ where
$0\leq r<k-1$.
Thus $d=nm-k+r+1$ and Inequality\eqnref{e-GriesMDS} for $i=0$ leads to
\[
  \frac{d}{q}\leq nm-d-\sum_{l=2}^{k-r-1}\Big\lceil\frac{d}{q^l}\Big\rceil
  \leq k-r-1-(k-r-2)=1,
\]
hence $q\geq d$.
\\
This covers all cases, since we always have $km\geq\delta$.
\end{proof}

The proof shows that in general the lower bounds on~$q$ are not tight since we have
estimated $\big\lceil\frac{d}{q^l}\big\rceil$ by~$1$ for $l\geq2$ in all cases.
For instance, if $(n-k+1)^2>d$, no \nkdqm-MDS code exists for $q=\frac{d}{n-k+1}$ and
$k=1$ or $km=\delta+1$.
But even if $\big\lceil\frac{d}{q^l}\big\rceil=1$ for all $l\geq2$ there might not
exist an \nkdq-MDS code where~$q$ attains the lower bound.
The obstacle is that for some $i\in\hat{\N}$ there might not exist an
$(n(m+i),k(m+i)-\delta)_q$-block code
with the appropriate distance as required by the proof of Theorem~\ref{T-CCB}.
Since these block codes have to produce a convolutional code in a very specific way,
they even have to have some additional structure.
We wish to illustrate this by the following example.

\begin{exa}\label{E-F3code}
Let $\nkd=(3,2,3)$. The generalized Singleton bound is $d:=S(3,2,3)=6$ and the
memory of a $(3,2,3)$-MDS code is $m=2$, see Theorem~\ref{T-MDSC}(b).
From Theorem~\ref{T-MDSfieldsize} we obtain $q\geq3$ for the field size.
Taking $q=3$ we have $\big\lceil\frac{d}{q^2}\big\rceil=1$ so that indeed the lower
bound for the field size cannot be improved.
The existence of a $(3,2,3;2)_3$-MDS code requires the existence of
$(3(2+i),1+2i)_3$-block codes with distance at least~$6$ for all $i\in\N_0$.
Such codes do indeed exist\footnote{
For small~$i$ these codes can be found in tables listing ternary codes.
For the general case we wish to thank H.-G.~Quebbemann who pointed out to us a
construction of such codes for sufficiently large~$i$ using direct products of
finitely many ``short'' MDS-codes over~$\F_{3^3}$ and mapping them into ternary codes.
}.
However, the block codes have to have some additional structure in order to be part
of a convolutional code.
To see this, let $G\in\F_3[z]^{2\times3}$ be a minimal generator matrix of the
desired convolutional code~$\cC$.
Write
\[
    G=\begin{bmatrix}g_1\\g_2\end{bmatrix}+z\begin{bmatrix}g_3\\g_4\end{bmatrix}
      +z^2\begin{bmatrix}g_5\\0\end{bmatrix}\text{ where }g_i\in\F_3^3.
\]
Recall from the proof of Theorem~\ref{T-CCB} that our arguments are based in particular
on the block code $\cC_1:=\{(u_1,u_2+u_3z)G\mid u_1,u_2,u_3\in\F_3\}$.
Comparing like powers of~$z$ one observes that this code is isomorphic to
\[
  \hat{\cC}_1=\im\begin{bmatrix}g_1&g_3&g_5\\g_2&g_4&0\\0&g_2&g_4\end{bmatrix}\subseteq\F_3^9.
\]
Using elementary row operations on the polynomial matrix~$G$ we may assume that the
entry of~$G$ at the position $(1,1)$ is a constant.
Furthermore, after rescaling the columns of~$G$ we may assume $g_4=(1,1,1)$.
Finally, due to non-catastrophicity, the entries of~$g_2$ are not all the same and
because of $\dist(\hat{\cC_1})=6$, all nonzero.
This gives us (up to block code equivalence) the two options
\[
  \im\begin{bmatrix}a_1&a_2&a_3&0&a_4&a_5&0&a_6&a_7\\1&1&2&1&1&1&0&0&0\\
                                0&0&0&1&1&2&1&1&1\end{bmatrix}
  \text{ or }\quad
  \im\begin{bmatrix}a_1&a_2&a_3&0&a_4&a_5&0&a_6&a_7\\1&2&2&1&1&1&0&0&0\\
                                0&0&0&1&2&2&1&1&1\end{bmatrix}
\]
for $\hat{\cC}_1$.
Going through some tedious calculations one can show that no such code in $\F_3^9$
with distance~$6$ exists.
Hence no $(3,2,3)_3$-MDS convolutional code exists.
\end{exa}

In the next section we will give examples of MDS codes over fields $\F_q$
where~$q$ attains the lower bound in all cases except for the case
$km=\delta+1$.

\section{Examples of Some Optimal Convolutional Codes}
\label{S-Exa}
\setcounter{equation}{0}

In this section we present some convolutional codes with distance reaching the
Griesmer bound.
To the best of our knowledge it was for most of the parameters, if not all,  not
known before whether such codes existed.

In the first column of the tables below the parameters of the given code are listed.
In the second column we give the Griesmer bound
$g:=G_q(n,k,\delta;m)$ for these parameters.
The third column gives a code reaching this bound.
In all examples the distance of the code has been
computed via a program.
In each case the code is given by a minimal generator matrix.
Thus, in particular all matrices given below are right invertible.
In the forth column we present the index of the first column distance that
reaches the free distance, cf.\eqnref{e-coldistfree}.
In the last column we indicate whether the code is a cyclic convolutional code in
the sense of Section~\ref{S-CCC}.
At the moment this additional structure is not important.
We only want to mention that cyclic convolutional codes do not exist for all sets of
parameters, in particular the length and the characteristic of the field have to be
coprime (just like for block codes).
Moreover, the shortest {\em binary\/} cyclic convolutional codes with
complexity $\delta>0$ have length $n=7$ or $n=15$.

The fields being used in the tables are $\F_2=\{0,1\},\,\F_4=\{0,1,\alpha,\alpha^2\}$
where $\alpha^2+\alpha+1=0$,
$\F_8=\{0,1,\beta,\ldots,\beta^6\}$ where $\beta^3+\beta+1=0$, and
$\F_{16}=\{0,1,\gamma,\ldots,\gamma^{14}\}$ where $\gamma^4+\gamma+1=0$.

The  generator matrix $\hat{G}_3$ of the $(15,4,12;3)_2$-code
in Table~I is given by
\[
  \hat{G}_3\T=\begin{bmatrix} 1+z^2&1+z+z^3&z+z^2&1+z+z^3 \\ 1+z+z^2&1+z+z^2+z^3&1+z+z^2+z^3&z \\
       1+z+z^3&1+z+z^2&1+z+z^2&1+z^2+z^3 \\
       z&1+z+z^3&1&1+z+z^2 \\ z&z^2&1+z&1+z^3 \\ z^2&z+z^3&z^3&1+z+z^2+z^3 \\ 1+z+z^3&z^2+z^3&z+z^2+z^3&z \\
       z^3&1+z+z^2&z+z^3&z^2 \\ z+z^2+z^3&z+z^2&1+z^3&z^2+z^3 \\ 1+z+z^2+z^3&z^2+z^3&z^2&1+z+z^2 \\
       1&1&z+z^2+z^3&z^2 \\ z^2+z^3&1+z&1&0 \\ 1+z&0&1+z^2+z^3&1+z^3 \\ z^2+z^3&1+z^2+z^3&z^3&1+z+z^3 \\
       1+z^2+z^3&z^3&1+z+z^2&z+z^2+z^3
      \end{bmatrix}.
\]

Some additional explanations and remarks will follow the tables.

\newpage
\begin{center}Table~I
\tabcolsep.4mm

\nopagebreak
{\footnotesize
\begin{turn}{90}
\mbox{}\hfill\begin{tabular}{|c|c|c|c|c|}\hline
 $\nkdqm$ & $g$ & code meeting the Griesmer bound & $d^c_i$ & \!cy\!
\\ \hline\\[-2.8ex]\hline\\[-3ex]\hline
 $(5,3,4;2)_2$ & $6$ &
   {\footnotesize$\begin{bmatrix}
      1+z^2 & 1+z & z & 1+z^2 & z+z^2 \\ 1+z & z & 1+z & 1 & z \\
      z & 1 & 1+z & 1+z & 1\end{bmatrix}$\quad (not even)}& $7$ &  \rule[-.6cm]{0cm}{1.4cm}
\\ \hline\\[-2.8ex]\hline\\[-3ex]\hline
 $(5,2,6;3)_2$ & $12$ &
 {\footnotesize$\begin{bmatrix}z^3+z^2+1\!&\!z^2+z\!&\!z^3+z+1\!&\!z^2+z\!&\!z^3+1\\
        z+1\!&\! z^3+z^2+1\!&\!z^3+z^2\!&\!z^3+z+1\!&\!z^2+z
  \end{bmatrix}$\quad (even)}& $10$ & \rule[-.4cm]{0cm}{1cm}
\\ \hline
 $(5,2,6;4)_2$ & $12$ &
 {\footnotesize$\begin{bmatrix}
  1+z^3+z^4\!&\!1+z+z^4\!&\!1+z^3\!&\!1+z^2+z^3\!&\!z+z^3+z^4\\
  1+z^2\!&\!1+z\!&\!z^2+z\!&\!z^2+z+1\!&\!z^2+z+1
  \end{bmatrix}$\quad (even)}& $10$ & \rule[-.4cm]{0cm}{1cm}
\\ \hline\\[-2.8ex]\hline\\[-3ex]\hline
 $(9,3,1;1)_8$ & $8^{*\bullet}$ &
 {\footnotesize$\begin{bmatrix}z+1&z+\beta&z&z+\beta^2&z+\beta^3&z+\beta^6&z+1&z&z+\beta\\
         1&\beta^2&\beta^5&\beta^6&\beta^6&\beta^5&\beta^2&1&0\\
         0&1&\beta^2&\beta^5&\beta^6&\beta^6&\beta^5&\beta^2&1
   \end{bmatrix}$} & $1$ & \rule[-.6cm]{0cm}{1.4cm}
\\ \hline\\[-2.8ex]\hline\\[-3ex]\hline
  $(3,2,2;1)_5$ & $5^{*\bullet}$ &
  {\footnotesize$\begin{bmatrix}2+3z&3z&4+4z\\4+2z&1+3z&2z\end{bmatrix}$}
  & $5$ &\rule[-.4cm]{0cm}{1cm}
\\ \hline\\[-2.8ex]\hline\\[-3ex]\hline
 $(7,3,3;1)_2$ & $8$ &
 {\footnotesize $G_1=\begin{bmatrix} 1&z&1+z&1+z&1&z&0\\z&1+z&0&1+z&1&1&z\\0&z&1&0&1+z&1+z&1+z
   \end{bmatrix}$\quad (even)} & $2$ & $\times$\rule[-.6cm]{0cm}{1.4cm}
\\ \hline
 $(7,3,6;2)_2$ & $12$ &
 {\footnotesize $G_2=\begin{bmatrix} 1+z^2&z+z^2&1+z&1+z&1+z^2&z&z^2\\z&1+z+z^2&0&1+z+z^2&1+z^2&1+z^2&z\\
         z^2&z+z^2&1+z^2&0&1+z&1+z+z^2&1+z
   \end{bmatrix}$\quad (even)} & $5$ & $\times$\rule[-.6cm]{0cm}{1.4cm}
\\ \hline
 $(7,3,9;3)_2$ & $16$ &
 {\footnotesize$\begin{bmatrix}
     1+z^2+z^3 & z+z^2 & 1+z+z^3 & 1+z & 1+z^2 & z+z^3 & z^2+z^3\\
     z & 1+z+z^2+z^3 & 0 & 1+z+z^2 & 1+z^2+z^3 & 1+z^2+z^3 & z+z^3 \\
     z^2+z^3 & z+z^2 & 1+z^2 & z^3 &1+z+z^3 & 1+z+z^2+z^3 & 1+z
   \end{bmatrix}$\quad (even?)} & $9$ & $\times$\rule[-.6cm]{0cm}{1.4cm}
\\ \hline
 $(7,3,12;4)_2$ & $20$ &
 {\footnotesize$\begin{bmatrix}
     1+z+z^3+z^4 \!\!&\!\! 1+z^3+z^4 \!\!&\!\! 1+z^2 \!\!&\!\! z+z^2+z^4 \!\!&\!\! 1+z^2+z^3 \!\!&\!\!
         z \!\!&\!\!  z+z^2+z^3+z^4\\
     z^2+z^3 \!\!&\!\! 1+z+z^2+z^4 \!\!&\!\! 1+z^4 \!\!&\!\! 1+z+z^2+z^3+z^4 \!\!&\!\! z \!\!&\!\!
         1+z+z^3+z^4 \!\!&\!\! z^2+z^3 \\
     z^2+z^4 \!\!&\!\! z \!\!&\!\! 1+z+z^3 \!\!&\!\! 1+z+z^2+z^4 \!\!&\!\! 1+z^2+z^3+z^4 \!\!&\!\!
         z^2+z^3+z^4 \!\!&\!\! 1+z+z^3
   \end{bmatrix}$\quad (doubly even?)}& $14$ & $\times$\rule[-.6cm]{0cm}{1.4cm}
\\ \hline\\[-2.8ex]\hline\\[-3ex]\hline
 $(15,4,4;1)_2$ & $16$ &
  {\footnotesize$\hat{G}_1=\left[\!\!\begin{array}{ccccccccccccccc}
   z& 0&z&1+z&0&0&1+z&1&0&1&z&1+z&1+z&1+z&1\\
  1&0&z&0& 1&0&z& 1+z&1+z&z&1&z& 1&1+z&1+z\\
  1&1& z& z&z&1+z&0&z&1&1+z&z& 1&0&1+z&1\\
  1+z&1+z &1 &z &0&z&1+z&0&0&1+z&1&0 &1 &z&1+z
  \end{array}\!\!\right]$\quad (even)} & $2$ & $\times$\rule[-.7cm]{0cm}{1.7cm}
\\ \hline
 $(15,4,8;2)_2$ & $24$ &
  {\footnotesize$\hat{G}_2\!\!=\!\!\left[\!\!\begin{array}{ccccccccccccccc}
       1+z^2\!\!&\!\!1+z+z^2\!\!&\!\!1+z\!\!&\!\!z\!\!\!&\!\!\!z\!\!&\!\!z^2\!\!&\!\! 1+z\!\!&\!\!
       0\!\!&\!\!z+z^2\!\!&\!\!
       1+z+z^2\!\!&\!\!1\!\!&\!\!z^2\!\!&\!\!1+z\!\!&\!\!z^2\!\!&\!\!\!1+z^2\\
       1+z\!\!&\!\!1+z+z^2\!\!&\!\!1+z+z^2\!\!\!&\!\!\!1+z\!\!&\!\!z^2\!\!&\!\! z\!\!&\!\!z^2\!\!&\!\!1+z+z^2\!\!&\!\!
       z+z^2\!\!&\!\!z^2\!\!&\!\!1\!\!&\!\!1+z\!\!&\!\!0\!\!&\!\!
       1+z^2\!\!\!&\!\!0\\
       z+z^2\!\!&\!\!1+z+z^2\!\!&\!\!1+z+z^2\!\!\!&\!\!1\!\!\!&\!\!1+z\!\!&\!\!0\!\!&\!\!z+z^2\!\!&\!\!z\!\!&\!\!
       1\!\!&\!\!z^2\!\!&\!\!z+z^2\!\!&\!\!1\!\!&\!\!1+z^2\!\!&\!\!0\!\!\!&\!\!1+z+z^2\\
       1+z\!\!&\!\!z\!\!&\!\!1+z^2\!\!\!&\!\!\!1+z+z^2\!\!\!&\!\!\!1\!\!&\!\!1+z+z^2 \!\!&\!\!z \!\!&\!\!z^2 \!\!&\!\!
       z^2\!\!&\!\!1+z+z^2 \!\!&\!\!z^2\!\!&\!\! 0 \!\!&\!\!1 \!\!&\!\!
       1+z \!\!\!&\!\!z+z^2
       \end{array}\!\!\!\right]$} & $5$ & $\times$\rule[-.8cm]{0cm}{1.8cm}
\\ & & \ \hfill(even?) & & \rule[-.1cm]{0cm}{0cm}\\ \hline
 $(15,4,12;3)_2$ & $32$ &
  {\footnotesize see $\hat{G}_3$ above,\quad (even?)} & & $\times$\rule[-.2cm]{0cm}{.7cm}
\\ \hline
\end{tabular}
\end{turn}
}
\end{center}
\vfill

\newpage
\begin{center}Table~II
\tabcolsep.4mm

\nopagebreak
\begin{turn}{90}
{\footnotesize
\begin{tabular}{|c|c|c|c|c|}\hline
 $\nkdqm$ & $g$ & code meeting the Griesmer bound & $d^c_i$ &\!cy\!
\\ \hline\\[-2.8ex]\hline\\[-3ex]\hline
 $(3,1,1;1)_4$ & $6^*$ &
  {\footnotesize$[\alpha+\alpha z,\ \alpha^2+\alpha z,\ 1+\alpha z]$}& $2^{**}$ &
  $\times$\rule[-.2cm]{0cm}{.7cm}
\\ \hline
 $(3,1,2;2)_4$ & $9^*$ &
  {\footnotesize$[\alpha+\alpha z+z^2,\ \alpha^2+\alpha z+\alpha^2 z^2,\
     1+\alpha z+\alpha z^2]$}& $5$ & $\times$\rule[-.2cm]{0cm}{.7cm}
\\ \hline
 $(3,1,3;3)_4$ & $12^{*\bullet}$ &
  {\footnotesize$[\alpha+\alpha z+z^2+\alpha^2 z^3,\ \alpha^2+\alpha z+\alpha^2 z^2+z^3,\
     1+\alpha z+\alpha z^2+\alpha z^3]$}& $7$ & $\times$\rule[-.2cm]{0cm}{.7cm}
\\ \hline
 $(3,1,4;4)_4$ & $14$ &
  {\footnotesize$[\alpha+\alpha z+z^2+\alpha^2 z^3+\alpha z^4,\
     \alpha^2+\alpha z+\alpha^2 z^2+z^3+\alpha z^4,\
     1+\alpha z+\alpha z^2+\alpha z^3+\alpha z^4]$}& $10$ & $\times$\rule[-.2cm]{0cm}{.7cm}
\\ \hline
 $(3,1,5;5)_4$ & $16$ &
  {\footnotesize$[\alpha+\alpha z+z^2+\alpha^2 z^3+\alpha z^4+\alpha z^5,\
     \alpha^2+\alpha z+\alpha^2 z^2+z^3+\alpha z^4+z^5,\
     1+\alpha z+\alpha z^2+\alpha z^3+\alpha z^4+\alpha^2 z^5]$}& $11$ & $\times$\rule[-.2cm]{0cm}{.7cm}
\\ \hline\\[-2.8ex]\hline\\[-3ex]\hline
 $(5,2,2;1)_4$ & $8$ &
   {\footnotesize$\begin{bmatrix}0& \alpha+z& \alpha^2+\alpha^2z& \alpha^2+\alpha^2z& \alpha+z\\
   \alpha+\alpha^2z& z& \alpha& \alpha^2+z& \alpha^2+\alpha^2z\end{bmatrix}$}&
   $2$& $\times$\rule[-.4cm]{0cm}{1.1cm}
\\ \hline
 $(5,2,4;2)_4$ & $12$ &
   {\footnotesize$\begin{bmatrix}
     0&\alpha+z+\alpha z^2& \alpha^2+\alpha^2z+\alpha^2z^2&\alpha^2+\alpha^2z+\alpha^2z^2&
     \alpha+z+\alpha z^2\\
     \alpha+\alpha^2z+\alpha z^2& z+\alpha^2z^2& \alpha+\alpha^2z^2&\alpha^2+z+\alpha z^2&
     \alpha^2+\alpha^2z\end{bmatrix}$}& $5$ & $\times$\rule[-.4cm]{0cm}{1.1cm}
\\ \hline
 $(5,2,6;3)_4$ & $16$ &
   {\footnotesize$\begin{bmatrix}
     0 &\alpha^2+\alpha^2 z +\alpha z^2+z^3 &1+\alpha z+\alpha^2z^2+\alpha^2z^3 &
     1+\alpha z+\alpha^2z^2+\alpha^2z^3 & \alpha^2+\alpha^2z+\alpha z^2+z^3 \\
     \alpha^2+\alpha z+\alpha z^2+\alpha^2z^3 & \alpha^2 z+\alpha^2z^2+\alpha^2z^3 &
     \alpha^2+\alpha^2z^2+z^3 & 1+\alpha^2 z+\alpha z^2 & 1+\alpha z+z^3
      \end{bmatrix}$}& $9$ & $\times$\rule[-.4cm]{0cm}{1.1cm}
\\ \hline\\[-2.8ex]\hline\\[-3ex]\hline
  $(3,2,2;1)_{16}$ &$5^*$&
  {\footnotesize$\begin{bmatrix}
      \gamma^5+\gamma^4 z & \gamma^3+\gamma^8 z & \gamma^9+\gamma^2 z\\
     \gamma^9+\gamma^{12} z & \gamma^5+\gamma^{14}z & \gamma^3+\gamma^3 z
    \end{bmatrix}$}
  & $3^{**}$ & $\times$\rule[-.4cm]{0cm}{1.1cm}
\\ \hline
  $(3,2,3;2)_{16}$ &$6^*$&
  {\footnotesize$\begin{bmatrix}
      \gamma+\gamma z+z^2 & \gamma^6+\gamma z+\gamma^{10}z^2 & \gamma^{11}+\gamma z+\gamma^{5}z^2\\
      1+z & \gamma^{10}+\gamma^{5}z & \gamma^5+\gamma^{10} z
    \end{bmatrix}$}
  & $5$ & $\times$\rule[-.4cm]{0cm}{1.1cm}
\\ \hline\\[-2.8ex]\hline\\[-3ex]\hline
 $(5,1,1;1)_{16}$ &$10^*$&
  {\footnotesize$[\gamma+\gamma z, \gamma^{13}+\gamma^{10} z, \gamma^{10}+\gamma^4 z,
     \gamma^7+\gamma^{13} z, \gamma^4+\gamma^{7}z]$}
  & $2^{**}$ & $\times$\rule[-.2cm]{0cm}{.7cm}
\\ \hline
 $(5,1,2;2)_{16}$ & $15^*$ &
  {\footnotesize$[\gamma+\gamma^4 z+\gamma z^2, \gamma^{7}+\gamma z+\gamma^{10}z^2,
     \gamma^{13}+\gamma^{13} z+\gamma^4z^2,\gamma^4+\gamma^{10}z+\gamma^{13}z^2,
     \gamma^{10}+\gamma^7z+\gamma^7z^2]$}
  & $3^{**}$ & $\times$\rule[-.2cm]{0cm}{.7cm}
\\ \hline
 $(5,1,3;3)_{16}$ & $20^*$ &
  {\footnotesize$[\gamma+z+\gamma^2 z^2+z^3, \gamma^{7}+\gamma^{12}z+\gamma^{11}z^2+\gamma^3z^3,
     \gamma^{13}+\gamma^{9}z+\gamma^5z^2+\gamma^6z^3,\gamma^4+\gamma^{6}z+\gamma^{14}z^2+\gamma^9z^3,
     \gamma^{10}+\gamma^3z+\gamma^8z^2+\gamma^{12}z^3]$}
  & $5$ & $\times$\rule[-.2cm]{0cm}{.7cm}
\\ \hline\\[-2.8ex]\hline\\[-3ex]\hline
 $(5,2,2;1)_{16}$ & $9^*$ &
   {\footnotesize$\begin{bmatrix}
      \gamma+\gamma z & \gamma^{13}+\gamma^{10} z & \gamma^{10}+\gamma^4 z & \gamma^7+\gamma^{13}z &
      \gamma^4+\gamma^7z\\
     1+\gamma^{5} z\ & \gamma^3+\gamma^{11}z\ & \gamma^6+\gamma^2 z & \gamma^9+\gamma^8 z &
     \gamma^{12}+\gamma^{14} z
    \end{bmatrix}$}
   & $2^{**}$ & $\times$\rule[-.4cm]{0cm}{1.1cm}
\\ \hline\\[-2.8ex]\hline\\[-3ex]\hline
  $(7,1,1;1)_{8}$ & $14^*$ &
   {\footnotesize$[\beta+\beta z,\, \beta^3+z,\, \beta^5+\beta^6z,\, 1+\beta^5z,\,
     \beta^2+\beta^4z,\, \beta^4+\beta^3z,\, \beta^6+\beta^2z]$}
   & $2^{**}$ & $\times$\rule[-.2cm]{0cm}{.7cm}
\\ \hline
  $(7,1,2;2)_{8}$ & $21^*$ &
   {\footnotesize$[\beta^2+\beta z+z^2, \beta^5+\beta^3z+\beta^6z^2, \beta+\beta^5z+\beta^5z^2,
       \beta^4+z+\beta^4z^2, 1+\beta^2z+\beta^3z^2, \beta^3+\beta^4z+\beta^2z^2,
       \beta^6+\beta^6z+\beta z^2]$}
  & $3^{**}$ & $\times$\rule[-.2cm]{0cm}{.7cm}
\\ \hline
  $(7,1,3;3)_{8}$ & $28^*$ &
   {\footnotesize$[1\!+\!\beta z\!+\!\beta^6z^2\!+\!z^3, 1\!+\!\beta^5z\!+\!\beta^5z^2\!+\!\beta^5z^3,
       1\!+\!\beta^2z\!+\!\beta^4z^2\!+\!\beta^3z^3,1\!+\!\beta^6z\!+\!\beta^3z^2\!+\!\beta z^3,
       1\!+\!\beta^3z\!+\!\beta^2z^2+\beta^6z^3, 1\!+\!z\!+\!\beta z^2+\beta^4z^3,
       1\!+\!\beta^4z\!+\!z^2\!+\!\beta^2z^3]$}
  & $5$ & $\times$\rule[-.2cm]{0cm}{.7cm}
\\ \hline
  $(7,2,3;2)_{8}$ & $14^*$ &
   {\footnotesize$\begin{bmatrix}
      1+z+\beta^4z^2 & \beta^4+\beta^5z+\beta^5z^2 & \beta+\beta^3z+\beta^6z^2 &
      \beta^5+\beta z+z^2 & \beta^2+\beta^6 z+\beta z^2 &
      \beta^6+\beta^4z+\beta^2z^2 & \beta^3+\beta^2z+\beta^3z^2\\
      \beta+\beta z & \beta^3+z & \beta^5+\beta^6z & 1+\beta^5z & \beta^2+\beta^4z &
      \beta^4+\beta^3z & \beta^6+\beta^2z
    \end{bmatrix}$}
  & $3$ & $\times$\rule[-.4cm]{0cm}{1.1cm}
\\ \hline
\end{tabular}
}
\end{turn}
\end{center}
\vfill

\newpage
\tabcolsep2mm
\begin{center}
Table~III
\\[1ex]
{\footnotesize
\begin{tabular}{|c|c|c|c|c|}\hline
 $\nkdqm$ & $g$ & code meeting the Griesmer bound & $d^c_i$ &\!cy\!
\\ \hline\\[-2.8ex]\hline\\[-3ex]\hline
  $(6,3,3;1)_{2}$ & $6$ & {\footnotesize columns 1,\,2,\,3,\,5,\,6,\,7 of $G_1$\quad (even)}
  & $3$ & \rule[-.2cm]{0cm}{.7cm}
\\ \hline
  $(6,3,6;2)_{2}$ & $10$ & {\footnotesize columns 1,\,2,\,4,\,5,\,6,\,7 of $G_2$\quad (even)}
  & $3$ & \rule[-.2cm]{0cm}{.7cm}
\\ \hline\\[-2.8ex]\hline\\[-3ex]\hline
  $(14,4,4;1)_{2}$ & $14$ & {\footnotesize columns 1 -- 14 of $\hat{G}_1$\quad (not even)}
  & $3$ & \rule[-.2cm]{0cm}{.7cm}
\\ \hline
  $(13,4,4;1)_{2}$ & $13$ & {\footnotesize columns 1,\,2,\,4 -- 14 of $\hat{G}_1$\quad (not even)}
  & $3$&  \rule[-.2cm]{0cm}{.7cm}
\\ \hline
  $(12,4,4;1)_{2}$ & $12$ & {\footnotesize columns 1,\,2,\,4 -- 12,\,14 of $\hat{G}_1$\quad (even)}
  & $3$ & \rule[-.2cm]{0cm}{.7cm}
\\ \hline
  $(10,4,4;1)_{2}$ & $10$ & {\footnotesize columns 1,\,2,\,4,\,6 -- 11,\,14 of $\hat{G}_1$\quad (even)}
  & $4$ & \rule[-.2cm]{0cm}{.7cm}
\\ \hline
  $(8,4,4;1)_{2}$ & $8$ & {\footnotesize columns 1,\,2,\,4,\,5,\,8,\,11,\,13,\,14 of $\hat{G}_1$\quad (not even)}
  & $4$ & \rule[-.2cm]{0cm}{.7cm}
\\ \hline\\[-2.8ex]\hline\\[-3ex]\hline
  $(14,4,8;2)_{2}$ & $22$ & {\footnotesize columns 2 -- 15 of $\hat{G}_2$\quad (even?)}
  & $6$ &\rule[-.2cm]{0cm}{.7cm}
\\ \hline
  $(13,4,8;2)_{2}$ & $20$ & {\footnotesize columns 1 -- 4,\,7 -- 15 of $\hat{G}_2$ \quad (even?)}
  & $6$ &\rule[-.2cm]{0cm}{.7cm}
\\ \hline
  $(12,4,8;2)_{2}$ & $18$ & {\footnotesize columns 1,\,2,\,4,\,7 -- 15 of $\hat{G}_2$\quad (not even)}
  & $6$ &\rule[-.2cm]{0cm}{.7cm}
\\ \hline
  $(10,4,8;2)_{2}$ & $16$ & {\footnotesize columns 1,\,2,\,4,\,5,\,7,\,8,\,10,\,11,\,13,\,14 of $\hat{G}_2$
  \quad (even?)}
  & $7$ &\rule[-.2cm]{0cm}{.7cm}
\\ \hline
  $(8,4,8;2)_{2}$ & $12$ & {\footnotesize columns 1,\,2,\,6,\,9,\,12 -- 15 of $\hat{G}_2$ (even?)}
  & $9$ & \rule[-.2cm]{0cm}{.7cm}
\\ \hline
\end{tabular}
}
\end{center}

It remains to explain some additional notation of the tables.
We also make some further comments illustrating the contents of the tables.

\begin{rem}\label{R-tables}
\begin{alphalist}
\item A $*$ attached to the bounds in the second column indicate that these numbers
      are identical to the generalized Singleton bound.
      Hence the corresponding codes are even MDS codes.
\item An additional supscript~${\bullet}$ attached to the bound~$g$
      indicates that the code is an MDS code where the field size reaches the
      lower bound of Theorem~\ref{T-MDSfieldsize}.
      This gives us examples for the three cases $k=1$, $km>\delta+1$, and
      $km=\delta$.
      We did not find an example of an \nkdq-MDS code
      where $km=\delta+1$ and $q=\frac{d}{n-k+1}$.
\item In~\cite[Prop.~2.3]{GRS03} it has been shown that the $j$th column distance of an
      \nkdq-code satisfies $d^c_j\leq(n-k)(j+1)+1$.
      From this it follows that the earliest column distance of an MDS code
      that can reach the free distance has index
      $M:=\big\lfloor\frac{\delta}{k}\big\rfloor+\big\lceil\frac{\delta}{n-k}\big\rceil$,
      see~\cite[Prop.~2.6]{GRS03}.
      In the same paper an MDS code is called strongly MDS if the $M$th column distance is
      equal to the free distance.
      We attached a~$^{**}$ to the index of the column distance in the second last column
      of the tables in order to indicate the strongly MDS codes.
      As far as we know no upper bound for the column distances is known that
      also takes the field sizes into account.
      However, using the estimate $d^c_j\leq(n-k)(j+1)+1$ one observes that the
      $(5,2,2;1)_4$- and the $(9,3,1;1)_8$-code are also optimal in the
      sense that no code with the same parameters exists where an earlier column distance
      reaches the free distance.
      We did not investigate whether any of the other codes is optimal in this sense.
\item We investigated the binary codes with respect of being even, that is, whether all codewords
      have even weight. This can be done by computing the weight
      distribution (see~\cite{McE98a} or~\cite[Sec.~3.10]{JoZi99}).
      Evenness of a code is indicated by an (even) attached to the generator matrix.
      Since the computation of the full weight distribution is very complex for larger
      complexity, we did not fully check the binary codes having complexity bigger than~$6$.
      In those cases we checked the weight of codewords associated with
      message words of small degree. In case this weight is always even we
      think there is strong evidence that the code is even and attached an
      (even?) to the generator matrix.
      In this sense there is also evidence that the $(7,3,12;4)_2$-code is doubly even,
      that is, all codewords have weight divisible by~$4$.
      Further investigation is necessary in order to understand whether (and why)
      all the binary cyclic convolutional codes of length~$7$ and~$15$ are even.
\item The second and third code of Table~I show that a code meeting the Griesmer bound
      need not have evenly distributed Forney indices.
      In other words, such a code need not be compact in the sense of Theorem~\ref{T-MDSC}(b).
      For both codes in Table~I the free distance is attained by the 10th column
      distance.
      Only the full weight distribution shows that the code with Forney
      indices $3,\,3$ is better than the code with indices $4,\,2$.
      The first one has weight distribution
      \[
        W_1(T)=10T^{12}+12T^{14}+71T^{16}+248T^{18}+873T^{20}+\ldots,
      \]
      saying that there are $10$ molecular codewords of weight~$12$ and $12$ molecular codewords
      of weight~$14$, etc. (for the definition of molecular codewords, see~\cite{McE98a};
      for weight distributions see also~\cite[Sec.~3.10]{JoZi99}).
      The weight distribution of the second code is
      \[
        W_2(T)=10T^{12}+27T^{14}+99T^{16}+350T^{18}+1280T^{20}+\ldots.
      \]
\item It is worth being mentioned that the codes with parameters
      $(7,3,3;1)_2,\,(7,3,6;2)_2$, and $(7,3,9;3)_2$ form a sequence in the sense
      that if one deletes $z^3$ (resp.~$z^2$) in the last (resp.~second) of the
      according generator matrices then one obtains the previous code.
      The same applies to the codes with parameters $(3,1,1;1)_4$, \ldots, $(3,1,5;5)_4$
      as well as to the $(5,2,2;1)_4$- and $(5,2,4;2)_4$-codes.
\item The codes with parameters $(7,3,3;1)_2,\,(7,3,6;2)_2,\,(15,4,4;1)_2$ and
      $(15,4,8;2)_2$ are extremely robust against puncturing in the sense of cutting
      columns of the according generator matrix
      (this is not puncturing in the sense of~\cite[Sec.~8]{McE98}).
      This way we do not only obtain right invertible matrices again, but even
      minimal matrices and, by doing this appropriately, codes reaching the Griesmer
      bound.
      We have cut one column of the codes of length~$7$ and up
      to~$7$ columns of the codes of length~$15$.
      The results are given in Table~III.
      The only cases where we did not get codes reaching the Griesmer bound
      are for $(11,4,4;1)_2$ and for $(9,4,8;2)_2$.
      We do not know if for these parameters there exist any codes at all that reach the
      bound.
      Since $G_2(9,4,4;1)=8=G_2(8,4,4;1)$ and
      $G_2(11,4,8;2)=16=G_2(10,4,8;2)$ we skipped in both cases the bigger length.
      Puncturing the code of length~$7$ and memory bigger
      than~$2$ did not result in a code meeting the Griesmer bound.
      We did not puncture the code of length~$15$ and memory~$3$.
\item Consider the $(8,4,4;1)_2$-code given in Table~III.
      There are other codes with exactly these parameters given in the literature.
      Indeed, in~\cite{JSW00} some (doubly-even self-dual) $(8,4,4;1)_2$-codes are presented.
      Our code is not even, which can
      easily be seen by writing down the generator matrix.
      We also computed the weight distribution and obtained
      \begin{align*}
        W(T)=&11T^8+28T^9+39T^{10}+101T^{11}+206T^{12}+565T^{13}+1374T^{14}+3033T^{15}\\
          &+7366T^{16}+16984T^{17}+40510T^{18}+95617T^{19}+22348T^{20}+\ldots,
      \end{align*}
      which is better than the weight distribution of the self-dual code given
      in~\cite[Eq.~(10)]{JSW00}.
\end{alphalist}
\end{rem}

\section{Cyclic Convolutional Codes}\label{S-CCC}
\setcounter{equation}{0}

The first two tables of the last section list plenty of optimal codes
that we have declared as cyclic.
Moreover, they gave rise to further sets of optimal codes as listed in
Table~III.
In this section we want to briefly describe the notion of cyclicity for
\CC{}s.
The first investigations in this direction have been made in the seventies
by Piret~\cite{Pi76} and Roos~\cite{Ro79}.
In both papers it has been shown (with different methods and in different
contexts) that cyclicity of \CC{}s must not be understood in
the usual sense, i.~e. invariance under the cyclic shift, if one wants to
go beyond the theory of cyclic block codes (see
Theorem~\ref{T-CCCclassic} below).
As a consequence, Piret suggested a more complex notion of cyclicity which then has been
further generalized by Roos.
In both papers some nontrivial examples of \CCC{}s in
this new sense are presented along with their distances.
All this indicates that the new notion of cyclicity seems to be the
appropriate one in the convolutional case.
Unfortunately, the papers~\cite{Pi76,Ro79} did not get much attention at
that time and the topic came to a halt.
Only recently it has been resumed in~\cite{GS02}.
Therein, an algebraic theory of \CCC{}s has been established which goes
well beyond the results of the seventies.
On the one hand it leads to a nice, yet nontrivial, generalization of the theory of
cyclic block codes, on the other hand it gives a very powerful toolbox for
constructing \CCC{}s.
We will now give a very brief description of these results and refer
to~\cite{GS02} for the details.

Just like for cyclic block codes we assume from now on that the length~$n$
and the field size~$q$ are coprime. Let $\F=\F_q$ be a field of size~$q$.
Recall that a block code $\cC\subseteq\F^n$ is called cyclic if it is
invariant under the cyclic shift, i.~e.
\begin{equation}\label{e-cs}
  (v_0,\ldots,v_{n-1})\in\cC\Longrightarrow
  (v_{n-1},v_0,\ldots,v_{n-2})\in\cC
\end{equation}
for all $(v_0,\ldots,v_{n-1})\in\F^n$.
It is well-known that this is the case if and only if~$\cC$ is an ideal in
the quotient ring
\begin{equation}\label{e-A}
     A:=\F[x]/_{\DS
     \ideal{x^n-1}}=\Big\{\sum_{i=0}^{n-1}f_ix^i\;\mod(x^n-1)\,\Big|\,f_0,\ldots,f_{n-1}\in\F\Big\},
\end{equation}
identified with $\F^n$ in the canonical way via
\[
  \p: \F^n\longrightarrow A,\quad
  (v_0,\ldots,v_{n-1})\longmapsto\sum_{i=0}^{n-1}v_ix^i.
\]
At this point it is important to recall that the cyclic shift in~$\F^n$ translates
into multiplication by~$x$ in~$A$, i.~e.
\begin{equation}\label{e-cx}
  \p(v_{n-1},v_0,\ldots,v_{n-2})=x\p(v_0,\ldots,v_{n-1})
\end{equation}
for all $(v_0,\ldots,v_{n-1})\in\F^n$.
Furthermore, it is well-known that each ideal $I\subseteq A$ is principal,
hence there exists some $g\in A$ such that $I=\ideal{g}$.
One can even choose~$g$ as a monic divisor of $x^n-1$, in which case it is
usually called the {\em generator polynomial\/} of the code
$\p^{-1}(I)\subseteq\F^n$.

It is our aim to extend this structure to the convolutional setting.
The most convenient way to do so is by using only the polynomial part
$\cC\cap\F[z]^n$ of the \CC\ $\cC\subseteq\Flaurent^n$.
Recall from\eqnref{e-cpolyunique} that this uniquely determines
the full code.
Hence imposing some additional structure on the polynomial part (that is,
on the generator matrix) will also impose some additional structure on
the full code.
In Remark~\ref{R-convcirc} below we will see from hindsight that one
can just as well proceed directly with the full code.
The polynomial part of a \CC\ is always a submodule of the free module
$\F[z]^n$.
Due to the right invertibility of the generator matrix not every submodule
of $\F[z]^n$ arises as polynomial part of a \CC.
It is easy to see~\cite[Prop.~2.2]{GS02} that we have

\begin{rem}\label{R-dirsumm}
A submodule $\cS\subseteq\F[z]^n$ is the polynomial part of some \CC\ if
and only if $\cS$ is a direct summand of $\F[z]^n$, i.e.
$\cS\oplus\cS'=\F[z]^n$ for some submodule $\cS'\subseteq\F[z]^n$.
\end{rem}

In order to extend the situation of cyclic block codes to the convolutional setting,
we have to replace the vector space~$\F^n$ by the free module $\F[z]^n$ and,
consequently, the ring~$A$ by the polynomial ring
\[
   A[z]:=\Big\{\sum_{j=0}^Nz^ja_j\,\Big|\, N\in\N_0,\,a_j\in A\Big\}
\]
over~$A$.
Then we can extend the map~$\p$ above coefficientwise to polynomials, thus
\begin{equation}\label{e-p}
  \p: \F[z]^n\longrightarrow A[z],\quad
  \sum_{j=0}^N z^jv_j\longmapsto \sum_{j=0}^N z^j\p(v_j),
\end{equation}
where, of course, $v_j\in\F^n$ and thus $\p(v_j)\in A$ for all~$j$.
This map is an isomorphism of $\F[z]$-modules.
Again, by construction the cyclic shift in $\F[z]^n$ corresponds to multiplication
by~$x$ in $A[z]$, that is, we have\eqnref{e-cx} for all
$(v_0,\ldots,v_{n-1})\in\F[z]^n$.
At this point it is quite natural to call a \CC\ $\cC\subseteq\Flaurent^n$ cyclic if it is
invariant under the cyclic shift, i.~e. if\eqnref{e-cs} holds true for all
$(v_0,\ldots,v_{n-1})\in\Flaurent^n$.
This, however, does not result in any codes other than block codes due to
the following result, see~\cite[Thm.~3.12]{Pi76} and \cite[Thm.~6]{Ro79}.
An elementary proof can be found at~\cite[Prop.~2.7]{GS02}.

\begin{theo}\label{T-CCCclassic}
Let $\cC\subseteq\Flaurent^n$ be an \nkd-convolutional code such that\eqnref{e-cs}
holds true for all $(v_0,\ldots,v_{n-1})\in\F[z]^n$.
Then $\delta=0$, hence~$\cC$ is a block code.
\end{theo}

This result has led Piret~\cite{Pi76} to suggest a different notion of
cyclicity for \CC{}s.
We will present this notion in the slightly more general version as it has been
introduced by Roos~\cite{Ro79}.

In order to do so notice that~$\F$ can be regarded as a subfield of the
ring~$A$ in a natural way. As a consequence,~$A$ is an $\F$-algebra, i.~e.,
a ring and a vector space over the field~$\F$ and the two structures are
compatible.
In the sequel the automorphisms of~$A$ with respect to this algebra structure
will play an important role. Therefore we define
\[
   \AutF(A):=\big\{\sigma:A\rightarrow
   A\,\big|\,\sigma|_{\F}=\text{id}_{\F},\,\sigma\text{ is bijective},\,
   \sigma(a\plusdot b)
   =\sigma(a)\plusdot\sigma(b)
   \text{ for all }a,\,b\in A
   \big\}.
\]
It is clear that each automorphism $\sigma\in\AutF(A)$ is uniquely
determined by the single value $\sigma(x)\in A$.
But not every choice for $\sigma(x)$ determines an automorphism on~$A$.
Since~$x$ generates the $\F$-algebra~$A$, the same has to be true
for~$\sigma(x)$ and, more precisely, we obtain for $a\in A$
\begin{equation}\label{e-sxa}
   \left.\begin{array}{l}
   \sigma(x)=a\text{ determines an }\\
   \text{automorphism on }A
   \end{array}\right\}
   \Longleftrightarrow
   \left\{\begin{array}{l}
   1,\,a,\ldots,a^{n-1}\text{ are linearly independent over }\F\\
   \text{and }a^n=1.
   \end{array}\right.
\end{equation}
Of course, $\sigma(x)=x$ determines the identity map on~$A$.
It should be mentioned that there is a better way to determine the automorphism
group of~$A$ by using the fact that the ring is direct product of fields.
This is explained in~\cite[Sec.~3]{GS02}.

The main idea of Piret was to impose a new ring structure on $A[z]$ and to call a
code cyclic if it is a left ideal with respect to that ring structure.
The new structure is non-commutative and based on an (arbitrarily chosen)
automorphism on~$A$.
In detail, this looks as follows.

\begin{defi}\label{D-CCC}
Let $\sigma\in\AutF(A)$.
\begin{arabiclist}
\item On the set $A[z]$ we define addition as usual and multiplication via
      \[
         \sum_{j=0}^Nz^j a_j\cdot\sum_{l=0}^M z^l b_l=
         \sum_{t=0}^{N+M}z^t\sum_{j+l=t}\sigma^l(a_j)b_l
         \text{ for all }N,\,M\in\N_0\text{ and }a_j,\,b_l\in A.
      \]
      This turns $A[z]$ into a non-commutative ring which is denoted by $\Azs$.
\item Consider the map $\p:\F[z]^n\!\rightarrow\!\Azs$ as in\eqnref{e-p}, where now
      the images $\p(v)\!=\!\sum_{j=0}^Nz^j\p(v_j)$ are regarded as elements of
      $\Azs$.
      A direct summand $\cS\subseteq\F[z]^n$ is said to be
      $\sigma$-{\em cyclic\/} if $\p(\cS)$ is a left ideal in $\Azs$.
\item A \CC\ $\cC\subseteq\Flaurent^n$ is said to be 
      $\sigma$-{\em cyclic\/} if $\cC\cap\F[z]^n$ is a \scy\ direct summand.
\end{arabiclist}
\end{defi}

A few comments are in order. First of all, notice that multiplication is
determined by the rule
\begin{equation}\label{e-az}
   az=z\sigma(a)\text{ for all }a\in A
\end{equation}
along with the rules of a (non-commutative) ring.
Hence, unless~$\sigma$ is the identity, the indeterminate~$z$ does not commute
with its coefficients.
Consequently, it becomes important to distinguish between left and right
coefficients of~$z$.
Of course, the coefficients can be moved to either side by applying the
rule\eqnref{e-az} since~$\sigma$ is invertible.
Multiplication inside~$A$ remains the same as before. Hence~$A$ is a
commutative subring of $\Azs$. Moreover, since
$\sigma|_{\F}=\text{id}_{\F}$, the classical polynomial ring $\F[z]$ is a
commutative subring of $\Azs$, too.
As a consequence, $\Azs$ is a left and right $\F[z]$-module and the
map~$\p:\F[z]^n\rightarrow\Azs$ is an isomorphism of left $\F[z]$-modules
(but not of right $\F[z]$-modules).
In the special case where $\sigma=\text{id}_{A}$ the ring $\Azs$ is
the classical commutative polynomial ring and we know from
Theorem~\ref{T-CCCclassic} that no \scy\ \CC{}s with nonzero complexity
exist.

\begin{exa}\label{E-binarylength7}
Let us consider the case where $\F=\F_2$ and $n=7$. Thus
$A=\F[x]/_{\ideal{x^7-1}}$.
Using\eqnref{e-sxa} one obtains~$18$ automorphisms, also listed
at~\cite[p.~680, Table~II]{Ro79}
(containing one typo: the last element of that table has to be
$x^2+x^3+x^4+x^5+x^6$ rather than $x+x^3+x^4+x^5+x^6$).
\\
Let us choose the automorphism $\sigma\in\AutF(A)$ defined by
$\sigma(x)=x^5$.
Furthermore, we consider the polynomial
\[
  g := 1+x^2+x^3+x^4+z(x+x^2+x^3+x^5)\in\Azs
\]
and denote by $\lideal{g}:=\{fg\mid f\in\Azs\}$ the left ideal generated
by~$g$ in $\Azs$.
Moreover, put $\cS:=\p^{-1}(\lideal{g})\subseteq\F[z]^7$.
We will show now that $\cS$ is a direct summand of $\F[z]^7$, hence
$\cS=\cC\cap\F[z]^7$ for some \CC\ $\cC\subseteq\Flaurent^7$, see
Remark~\ref{R-dirsumm}.
In order to do so we first notice that
\[
   \lideal{g}=\spann_{\F[z]}\big\{g,\,xg,\ldots,x^6g\big\}
\]
and therefore
\[
  \cS=\big\{u M\,\big|\, u\in\F[z]^7\big\}
  \text{ where }
       M=\begin{bmatrix}\p^{-1}(g)\\\p^{-1}(xg)\\\vdots\\\p^{-1}(x^6g)
       \end{bmatrix}.
\]
Thus we have to compute $x^ig$ for $i=1,\ldots,6$.
Using the multiplication rule in\eqnref{e-az} we obtain
\begin{align*}
  xg&=x+x^3+x^4+x^5+z(1+x+x^3+x^6),
  \\
  x^2g&=x^2+x^4+x^5+x^6+z(x+x^4+x^5+x^6),
  \\
   x^3g&=1+x^3+x^5+x^6+z(x^2+x^3+x^4+x^6)\\
       &=g+x^2g.
\end{align*}
Since $x^3g$ is in the $\F$-span of the previous elements, we obtain
$\lideal{g}=\spann_{\F[z]}\big\{g,xg,x^2g\big\}$ and,
since~$\p$ is an isomorphism,
\[
  \cS=\big\{u G\,\big|\, u\in\F[z]^3\big\},
\]
where
\[
  G=\begin{bmatrix}\p^{-1}(g)\\\p^{-1}(xg)\\\p^{-1}(x^2g)
    \end{bmatrix}
   =\begin{bmatrix} 1&z&1+z&1+z&1&z&0\\z&1+z&0&1+z&1&1&z\\0&z&1&0&1+z&1+z&1+z
    \end{bmatrix}.
\]
One can easily check that the matrix~$G$ is right invertible.
Hence $\cS$ is indeed a direct summand of $\F[z]^7$ and thus we have obtained a
\scy\ \CC\ $\cC=\im G\subseteq\Flaurent^7$.
This is exactly the $(7,3,3;1)_2$-code given in Table~I of the last
section.
\end{exa}

The other \CCC{}s in Tables~I and~II are obtained in a similar way.
Since the underlying automorphism cannot easily be read off from the generator
matrix of a \CCC\, we will, for sake of completeness, present them explicitly
in the following table.
All those codes come from principal left ideals in $\Azs$ and, except for the codes with
parameters $(3,2,3;2)_{16},\,(5,2,2;1)_{16},\,(7,2,3;2)_8$, the generator polynomial
can be recovered from the given data by applying the map~$\p$ to the first row of
the respective generator matrix.
The generator matrices of the remaining three codes are built in a slightly different
way.
In those cases each row of the given matrix generates a 1-dimensional cyclic code
and thus each of those three codes is the direct sum of two 1-dimensional cyclic
codes.
In each case a generator polynomial of the associated principal left ideal is obtained by
applying~$\p$ to the sum of the two rows of the respective generator matrix.

\begin{samepage}
\begin{center}
Table~IV
\\[1ex]
{\footnotesize
\begin{tabular}{|c|c|}\hline
 $\nkdqm$-code of Tables~I~and~II & automorphism given by \rule[-.2cm]{0cm}{.7cm}
\\ \hline\\[-2.8ex]\hline\\[-3ex]\hline
  $(7,3,3m;m)_2,\,m=1,\ldots,4$ & $\sigma(x)=x^5$  \rule[-.2cm]{0cm}{.7cm}
\\ \hline
  $(15,4,4;1)_2$ & $\sigma(x)=x+x^7+x^{10}$  \rule[-.2cm]{0cm}{.7cm}
\\ \hline
  $(15,4,4m;m)_2,\,m=2,3$ & $\sigma(x)=x^3+x^5+x^7+x^{10}+x^{12}+x^{13}+x^{14}$  \rule[-.2cm]{0cm}{.7cm}
\\ \hline
  $(3,1,\delta;\delta)_4,\, \delta=1,\ldots,5$ & $\sigma(x)=\alpha^2x$  \rule[-.2cm]{0cm}{.7cm}
\\ \hline
  $(5,2,2m;m)_4,\,m=1,2,3$ & $\sigma(x)=x^2$  \rule[-.2cm]{0cm}{.7cm}
\\ \hline
  $(3,2,2;1)_{16}$ and $(3,2,3;2)_{16}$ & $\sigma(x)=\gamma^{10}x$  \rule[-.2cm]{0cm}{.7cm}
\\ \hline
  $(5,1,\delta;\delta)_{16},\,\delta=1,2,3$ and $(5,2,2;1)_{16}$ & $\sigma(x)=x^3$  \rule[-.2cm]{0cm}{.7cm}
\\ \hline
  $(7,1,\delta;\delta)_{8},\,\delta=1,2$ and $(7,2,3;2)_{8}$ & $\sigma(x)=x^5$  \rule[-.2cm]{0cm}{.7cm}
\\ \hline
  $(7,1,3;3)_{8}$ & $\sigma(x)=\beta x+\beta x^2+\beta^3x^3+\beta^3x^4+\beta^3x^5+\beta^2x^6$  \rule[-.2cm]{0cm}{.7cm}
\\ \hline
\end{tabular}
}
\end{center}
\end{samepage}

The fact that all the \CCC{}s above come from principal left ideals in $\Azs$ is not a restriction since
we have the following important result.

\begin{theo}\label{T-pli}
Let $\sigma\in\AutF(A)$.
If $\cS\subseteq\F[z]^n$ is a \scy\ direct summand, then $\p(\cS)$ is a
principal left ideal of $\Azs$, that is, there exists some polynomial
$g\in\Azs$ such that $\p(\cS)=\lideal{g}$.
We call~$g$ a generator polynomial of both $\cS$ and the \scy\ \CC\
$\cC\subseteq\Flaurent^n$ determined by $\cS$, see Remark~\ref{R-dirsumm}
and\eqnref{e-cpolyunique}.
\end{theo}

The generator polynomial of a \scy\ \CC\ can be translated into vector notation and
leads to a generalized circulant matrix.
This looks as follows.
Let $\cS\subseteq\F[z]^n$ be a \scy\ direct summand and let $\p(\cS)=\lideal{g}$.
Define
\[
  \Msigma(g)=\begin{bmatrix}
    \p^{-1}(g)\\ \p^{-1}(xg)\\ \vdots\\
    \p^{-1}(x^{n-1}g)\end{bmatrix}\in\F[z]^{n\times n}.
\]
Then it is easy to see that $\p\big(u\Msigma(g)\big)=\p(u)g$
for all $u\in\F[z]^n$ (see~\cite[Prop.~6.8(b)]{GS02}) and therefore,
$\cS=\big\{u\Msigma(g)\,\big|\,u\in\F[z]^n\big\}$.
We call $\Msigma(g)$ the $\sigma$-circulant associated with~$g$.

\begin{rem}\label{R-convcirc}
Using the identities above we can now easily see that \scy\ structure can
also be considered without restricting to the polynomial part.
Just like the polynomial ring $A[z]$ we can turn the set $A(\!(z)\!)$ of formal
Laurent series over~$A$ into a non-commutative ring by
defining addition as usual and multiplication via\eqnref{e-az}.
We will denote the ring obtained this way by $\Azzs$.
Furthermore, we can extend the map $\p$ to Laurent series in the canonical way,
see also\eqnref{e-p}.
Then one can easily show that just like in the polynomial case
\[
   \p\big(u\Msigma(g)\big)=\p(u)g\text{ for all }u\in\Flaurent^n
\]
for each $g\in\Azs$.
Using the fact that a code $\cC\subseteq\Flaurent^n$ is uniquely determined
by its polynomial part (see\eqnref{e-cpolyunique}), and that the latter
is a principal left ideal in $\Azs$ due to Theorem~\ref{T-pli},
one can now derive the equivalence
\[
  \cC\subseteq\Flaurent^n\text{ is $\sigma$-cyclic }
  \Longleftrightarrow
  \p(\cC)\text{ is a left ideal in }\Azzs.
\]
Moreover, if $\cC$ is $\sigma$-cyclic, a generator polynomial of the ideal
$\p(\cC\cap\F[z]^n)$ in $\Azs$ is also a principal generator of the
ideal $\p(\cC)$ in $\Azzs$.
This justifies to call~$g$ a generator polynomial of
the full code~$\cC$ as we did in Theorem~\ref{T-pli}.
\end{rem}

At this point the question arises as to how a (right invertible)
generator matrix can be obtained from the $\sigma$-circulant $\Msigma(g)$.
Notice that in Example~\ref{E-binarylength7} the generator matrix of the code is
simply given by the first three rows of the circulant.
This is indeed in general the case, but requires a careful choice of the generator
polynomial~$g$ of the code.
Recall that, due to zero divisors in $\Azs$, the generators of a principal left
ideal, are highly non unique.
The careful choice of the generator polynomial is based on a Gr\"obner basis
theory that can be established in the non-commutative polynomial ring $\Azs$.
This is a type of reduction procedure resulting in unique generating sets of left
ideals which in turn produce very powerful $\sigma$-circulants.
The details of this theory goes beyond the scope of this paper and we refer the
reader to~\cite{GS02} for the details, in particular to
\cite[Thm.~7.8,~Thm.~7.18]{GS02}.
Therein it has been shown that a reduced generator polynomial also reflects the
parameters of the code, i.~e., the dimension and the complexity, and even leads to a
minimal generator matrix through $\sigma$-circulants.
Only with these results it becomes clear that \CCC{}s can have only very specific parameters
(length, dimension, and complexity) depending on the chosen field $\F_q$.
Furthermore, the notions of parity check polynomial and associated parity check
matrix have been discussed in detail in~\cite{GS02}, leading to a generalization of
the block code situation.

As for the cyclic codes of the last section we only would like to mention that their
generator polynomials obtained as explained right before Table~IV are all reduced in
the sense above.

So far we do not have any estimates for the distance of a \CCC\ in terms of its
(reduced) generator polynomial and the chosen automorphism.
The examples given in the last section have been found simply by trying some
promising reduced generator polynomials (using the algebraic theory of~\cite{GS02}).
Except for the puncturing in Table~III we did not perform a systematic search for optimal
codes.

\section*{Conclusion}
In this paper we gave many examples of \CCC{}s that all reach the Griesmer bound.
The examples indicate that this class of \CC{}s promises to contain many excellent
codes and therefore deserves further investigation.
As one of the next steps the relation between the (reduced) generator
polynomial and the automorphism on the one hand and the distance on
the other hand should be investigated in detail.

\bibliographystyle{abbrv}
\bibliography{literatureAK,literatureLZ}
\end{document}